\documentclass[12pt]{article}
\usepackage{amssymb,amsmath,amsthm,amsfonts,enumerate, bbm, mathdots}

\usepackage{latexsym}
\usepackage{amscd}
\usepackage{stmaryrd}
\usepackage{color}
\usepackage[all]{xypic}
\usepackage{epsfig}
\usepackage{graphics}
\usepackage{ifthen}
\usepackage{varioref}
\usepackage{rotating}
\usepackage{extarrows}
\usepackage{cite}
\usepackage{mathrsfs}

\numberwithin{equation}{section}

\theoremstyle{plain}
\newtheorem{thm}{Theorem}[section]
\newtheorem{cor}[thm]{Corollary}

\newtheorem{prop}[thm]{Proposition}



\definecolor{darkgreen}{rgb}{0.0625,0.64,0.0625}
\usepackage[
            pdfstartview=FitH,
            CJKbookmarks=true,
            bookmarksnumbered=true,
            bookmarksopen=true,
            colorlinks,
            pdfborder=001,
            linkcolor=blue,
            anchorcolor=green,
            citecolor=red
            ]{hyperref}

\def\proof{\noindent {\bf Proof.\;}}

\allowdisplaybreaks

\begin{document}
\title{Ohno-Zagier type relation for multiple $t$-values\thanks{This work was supported by the Fundamental Research Funds for the Central Universities (grant number 22120210552).}}
\date{\small ~ \qquad\qquad School of Mathematical Sciences, Tongji University \newline No. 1239 Siping Road,
Shanghai 200092, China}

\author{Zhonghua Li\thanks{E-mail address: zhonghua\_li@tongji.edu.cn} ~and  Yutong Song\thanks{E-mail address: 1853404@tongji.edu.cn}}
\maketitle
\begin{abstract}
We study the Ohno-Zagier type relation for multiple $t$-values and multiple $t$-star values. We represent the generating function of sums of multiple $t$-(star) values with fixed weight, depth and height in terms of the generalized hypergeometric function $\,_3F_2$. As applications, we get a formula for the generating function of sums of multiple $t$-(star) values of maximal height and a weighted sum formula for sums of multiple $t$-(star) values with fixed weight and depth.
\end{abstract}

{\small

}

\section{Introduction}

A finite sequence $\mathbf{k}=(k_1,\ldots,k_n)$ of positive integers is called an index. The weight, depth and height of the index $\mathbf{k}$ are defined  respectively by $k_1+\cdots+k_n$, $n$ and $|\{j\mid 1\leq j\leq n, k_j\geq 2\}|$. If $k_1>1$, the index $\mathbf{k}$ is called admissible. For an admissible index $\mathbf{k}=(k_1,\ldots,k_n)$, the multiple zeta value $\zeta(\mathbf{k})$ and the multiple zeta-star value $\zeta^{\star}(\mathbf{k})$ are defined respectively by
\begin{align*}
&\zeta(\mathbf{k})=\zeta(k_1,\ldots,k_n)=\sum\limits_{m_1>\cdots>m_n>0}\frac{1}{m_1^{k_1}\cdots m_n^{k_n}},\\
&\zeta^{\star}(\mathbf{k})=\zeta^{\star}(k_1,\ldots,k_n)=\sum\limits_{m_1\geq\cdots\geq m_n>0}\frac{1}{m_1^{k_1}\cdots m_n^{k_n}}.
\end{align*}
The systematic study of multiple zeta values was carried out by  M. E. Hoffman \cite{Hoffman92} and D. Zagier \cite{Zagier}. More results of multiple zeta values can be gotten from the book \cite{Zhao} of J. Zhao.

In this paper, we focus on the Ohno-Zagier type relations. For nonnegative integers $k,n,s$, denote by $I_0(k,n,s)$ the set of admissible indices of weight $k$, depth $n$ and height $s$. It is easy to see that $I_0(k,n,s)$ is nonempty if and only if $k\geq n+s$ and $n\geq s\geq 1$. Using the Gaussian hypergeometric function, Y. Ohno and D. Zagier proved in \cite{Ohno-Zagier} that
\begin{align*}
&\sum\limits_{k\geq n+s,n\geq s\geq 1}\left(\sum\limits_{\mathbf{k}\in I_0(k,n,s)}\zeta(\mathbf{k})\right)u^{k-n-s}v^{n-s}w^{s-1}\\
=&\frac{1}{uv-w}\left\{1-\exp\left(\sum\limits_{n=2}^\infty\frac{\zeta(n)}{n}(u^n+v^n-\alpha^n-\beta^n)\right)\right\},
\end{align*}
where $\alpha$ and $\beta$ are determined by $\alpha+\beta=u+v$ and $\alpha\beta=w$.

After that, similar studies were carried out on various generalizations of multiple zeta values. For example, in \cite{AKO}, a similar formula holds for the sums of multiple zeta-star values. While this time we do not meet the Gaussian hypergeometric function but the generalized hypergeometric function $\,_3F_2$. In fact, T. Aoki, Y. Kombu and Y. Ohno showed that
\begin{align*}
&\sum\limits_{k\geq n+s,n\geq s\geq 1}\left(\sum\limits_{\mathbf{k}\in I_0(k,n,s)}\zeta^{\star}(\mathbf{k})\right)u^{k-n-s}v^{n-s}w^{2s-2}\\
=&\frac{1}{(1-v)(1-\beta)}\,_3F_2\left({1-\beta,1-\beta+u,1\atop 2-v,2-\beta};1\right)
\end{align*}
with $\alpha$ and $\beta$ determined by $\alpha+\beta=u+v$ and $\alpha\beta=uv-w^2$. Here for a positive integer $m$ and complex numbers $a_1,\ldots,a_{m+1},b_1,\ldots,b_m$ with $b_1,\ldots,b_m\neq 0,-1,-2,\ldots$, the generalized hypergeometric function $\,_{m+1}F_m$ is defined as
\begin{align*}
\,_{m+1}F_m\left({a_1,\ldots,a_{m+1}\atop b_1,\ldots,b_m};z\right)=&\sum\limits_{n=0}^\infty\frac{(a_1)_n\cdots(a_{m+1})_n}{n!(b_1)_n\cdots(b_m)_n}z^n,
\end{align*}
with the Pochhammer symbol $(a)_n$ defined by
$$(a)_n=\frac{\Gamma(a+n)}{\Gamma(a)}=\begin{cases}
1, & \text{if\;} n=1,\\
a(a+1)\cdots (a+n-1), & \text{if\;} n>1.
\end{cases}$$
It is known that this formal power series converges absolutely for $|z|<1$, and it also converges absolutely for $|z|=1$ if $\Re(\sum b_i-\sum a_i)>0$. If $m=1$, we get the Gaussian hypergeometric function.

In \cite{Takeyama}, Y. Takeyama studied the Ohno-Zagier type relation for a level two variant of multiple zeta values, which is called multiple $T$-values and was introduced by M. Kaneko and H. Tsumura \cite{Kaneko-Tsumura}. As a consequence, a weighted sum formula of multiple $T$-values with fixed weight and depth was given in \cite{Takeyama}.  In this paper, we consider another level two variant of multiple zeta values, which is called multiple $t$-values and was introduced by M. E. Hoffman in  \cite{Hoffman}.

For an admissible index $\mathbf{k}=(k_1,\ldots,k_n)$, the multiple $t$-value $t(\mathbf{k})$ and the multiple $t$-star value $t^{\star}(\mathbf{k})$ are defined respectively by
\begin{align*}
&t(\mathbf{k})=t(k_1,\ldots,k_n)=\sum\limits_{m_1>\cdots>m_n>0\atop m_i:\text{odd}}\frac{1}{m_1^{k_1}\cdots m_n^{k_n}},\\
&t^{\star}(\mathbf{k})=t^{\star}(k_1,\ldots,k_n)=\sum\limits_{m_1\geq \cdots\geq m_n>0\atop m_i:\text{odd}}\frac{1}{m_1^{k_1}\cdots m_n^{k_n}}.
\end{align*}
It is easy to obtain the following iterated integral representations:
\begin{align*}
&t(\mathbf{k})=\int_0^1\left(\frac{dt}{t}\right)^{k_1-1}\frac{tdt}{1-t^2}\cdots \left(\frac{dt}{t}\right)^{k_{n-1}-1}\frac{tdt}{1-t^2}\left(\frac{dt}{t}\right)^{k_n-1}\frac{dt}{1-t^2},\\
&t^{\star}(\mathbf{k})=\int_0^1\left(\frac{dt}{t}\right)^{k_1-1}\frac{dt}{t(1-t^2)}\cdots\left(\frac{dt}{t}\right)^{k_{n-1}-1}\frac{dt}{t(1-t^2)}\left(\frac{dt}{t}\right)^{k_n-1}\frac{dt}{1-t^2},
\end{align*}
where $\left(\frac{dt}{t}\right)^{k_i-1}=\underbrace{\frac{dt}{t}\cdots\frac{dt}{t}}_{k_i-1}$, and for one forms $\omega_i(t)=f_i(t)dt$, $i=1,2,\ldots,k$, we define the iterated integral
$$\int_0^z\omega_1(t)\cdots\omega_k(t)=\int\limits_{z>t_1>\cdots>t_k>0}f_1(t_1)\cdots f_k(t_k)dt_1\cdots dt_k.$$

We want to study the sum of multiple $t$-(star) values with fixed weight, depth and height.  For nonnegative integers $k,n,s$, define
$$G_0(k,n,s)=\sum\limits_{\mathbf{k}\in I_0(k,n,s)}t(\mathbf{k}),\quad G_0^{\star}(k,n,s)=\sum\limits_{\mathbf{k}\in I_0(k,n,s)}t^{\star}(\mathbf{k}).$$
Then we obtain Ohno-Zagier type relations which  represent the generating functions of $G_0(k,n,s)$ and $G_0^{\star}(k,n,s)$ by the generalized hypergeometric function $\,_3F_2$.

\begin{thm}\label{Thm:G-0-kns}
For formal variables $u,v,w$, we have
$$\sum\limits_{k\geq n+s,n\geq s\geq 1}G_0(k,n,s)u^{k-n-s}v^{n-s}w^{s-1}=\frac{1}{1-u}\,_3F_2\left({\alpha,\beta,1\atop \frac{3-u}{2},\frac{3}{2}};1\right),$$
where $\alpha$ and $\beta$ are determined by $\alpha+\beta=1-\frac{u}{2}+\frac{v}{2}$ and $\alpha\beta=\frac{1}{4}(1-u+v-uv+w)$.
\end{thm}

\begin{thm}\label{Thm:G-0-kns-Star}
For formal variables $u,v,w$, we have
$$\sum\limits_{k\geq n+s,n\geq s\geq 1}G_0^{\star}(k,n,s)u^{k-n-s}v^{n-s}w^{s-1}=\frac{1}{1-u-v+uv-w}\,_3F_2\left({\frac{1-u}{2},\frac{1}{2},1\atop \alpha^{\star},\beta^{\star}};1\right),$$
where $\alpha^{\star}$ and $\beta^{\star}$ are determined by $\alpha^{\star}+\beta^{\star}=3-\frac{u}{2}-\frac{v}{2}$ and $\alpha\beta=\frac{1}{4}(9-3u-3v+uv-w)$.
\end{thm}

From the above theorems and using some summation formulas of $\,_3F_2$, we obtain some corollaries in Section  \ref{Sec:Application}. For example, we give a formula for the generating function of sums of multiple $t$-(star) values of maximal height and a weighted sum formula of sums of multiple $t$-(star) values with fixed weight and depth. Finally, we prove Theorem \ref{Thm:G-0-kns} and Theorem \ref{Thm:G-0-kns-Star} in Section \ref{Sec:Proof}.


\section{Applications}\label{Sec:Application}

\subsection{Sum of height one}

To save space, we denote a sequence of $k$ repeated $n$ times by $\{k\}^n$.

Setting $w=0$ in Theorem \ref{Thm:G-0-kns}, we find
$$\alpha,\beta=\frac{1-u}{2},\frac{1+v}{2}.$$
Hence we get the following result, which gives the generating function of height one multiple $t$-values.

\begin{cor}[{\cite[Theorem 5.1]{Hoffman}}]
We have
$$\sum\limits_{k\geq n+1,n\geq 1}t(k-n+1,\{1\}^{n-1})u^{k-n-1}v^{n-1}=\frac{1}{1-u}\,_3F_2\left({\frac{1-u}{2},\frac{1+v}{2},1\atop \frac{3-u}{2},\frac{3}{2}};1\right).$$
\end{cor}

Similarly, setting $w=0$ in Theorem \ref{Thm:G-0-kns-Star}, we have
$$\alpha^{\star},\beta^{\star}=\frac{3-u}{2},\frac{3-v}{2}.$$
Therefore, we obtain the generating function of height one multiple $t$-star values.

\begin{cor}
We have
$$\sum\limits_{k\geq n+1,n\geq 1}t^{\star}(k-n+1,\{1\}^{n-1})u^{k-n-1}v^{n-1}=\frac{1}{(1-u)(1-v)}\,_3F_2\left({\frac{1-u}{2},\frac{1}{2},1\atop \frac{3-u}{2},\frac{3-v}{2}};1\right).$$
\end{cor}

Note that using \cite[Lemma 5.2]{Hoffman}, for any integer $m\geq 2$, we have
\begin{align*}
\sum\limits_{n=1}^\infty t(m,\{1\}^{n-1})v^{n-1}=\,_{m+1}F_{m}\left({\frac{1+v}{2},1,\left\{\frac{1}{2}\right\}^{m-1}\atop \left\{\frac{3}{2}\right\}^m};1\right),\\
\sum\limits_{n=1}^\infty t^{\star}(m,\{1\}^{n-1})v^{n-1}=\frac{1}{1-v}\,_{m+1}F_{m}\left({1,\left\{\frac{1}{2}\right\}^{m}\atop \frac{3-v}{2},\left\{\frac{3}{2}\right\}^{m-1}};1\right).
\end{align*}

\subsection{Sum of maximal height}

Setting $v=0$ in Theorem \ref{Thm:G-0-kns}, we get the generating function of sums of multiple $t$-values of maximal height. By the symmetric sum formula \cite[Theorem 3.2]{Hoffman}, the sum of multiple $t$-values with fixed weight, depth and maximal height can be represented by $t$-values. Here we give a closed formula for the generating function of the sums of maximal height.

\begin{cor}
For formal variables $u,w$, we have
\begin{align}
1+\sum\limits_{k\geq 2n,n\geq 1}G_0(k,n,n)u^{k-2n}w^n=\exp\left\{\sum\limits_{n=2}^\infty\frac{t(n)}{n}(u^n-x^n-y^n)\right\},
\label{Eq:MaxHeight}
\end{align}
where $x$ and $y$ are determined by $x+y=u$ and $xy=w$.
\end{cor}

\proof Setting $v=0$ in Theorem \ref{Thm:G-0-kns}, we get $\alpha+\beta=1-\frac{u}{2}$ and $\alpha\beta=\frac{1}{4}(1-u+w)$. Let $x=1-2\alpha$ and $y=1-2\beta$, then $x+y=u$ and $xy=w$. Using the summation formula \cite[7.4.4.28]{PBA}
\begin{align}
&\,_3F_2\left({a,b,1\atop c,2+a+b-c};1\right)\nonumber\\
=&\frac{1+a+b-c}{(1+a-c)(1+b-c)}\left(1-c+\frac{\Gamma(c)\Gamma(1+a+b-c)}{\Gamma(a)\Gamma(b)}\right)
\label{Eq:3F2-1}
\end{align}
with $a=\alpha$, $b=\beta$ and $c=\frac{3}{2}$, we find that
\begin{align*}
&\frac{1}{1-u}\,_3F_2\left({\alpha,\beta,1\atop \frac{3-u}{2},\frac{3}{2}};1\right)
=\frac{1}{1-u}\frac{\frac{1-u}{2}}{\left(\alpha-\frac{1}{2}\right)\left(\beta-\frac{1}{2}\right)}\left(-\frac{1}{2}+\frac{\Gamma\left(\frac{3}{2}\right)\Gamma\left(\frac{1-u}{2}\right)}{\Gamma(\alpha)\Gamma(\beta)}\right)\\
=&-\frac{1}{w}+\frac{\sqrt{\pi}}{w}\frac{\Gamma\left(\frac{1-u}{2}\right)}{\Gamma\left(\frac{1-x}{2}\right)\Gamma\left(\frac{1-y}{2}\right)}.
\end{align*}
Using the duplication formula
$$\Gamma\left(\frac{1}{2}-\frac{z}{2}\right)=\frac{\sqrt{\pi}2^{z}\Gamma(1-z)}{\Gamma\left(1-\frac{z}{2}\right)},$$
and the expansion
\begin{align}
\Gamma(1-z)=\exp\left(\gamma z+\sum\limits_{n=2}^\infty\frac{\zeta(n)}{n}z^n\right),
\label{Eq:Gamma}
\end{align}
where $\gamma$ is Euler's constant, we get
$$\Gamma\left(\frac{1-z}{2}\right)=\sqrt{\pi}2^z\exp\left\{\frac{\gamma z}{2}+\sum\limits_{n=2}^\infty \frac{\zeta(n)}{n}(1-2^{-n})z^n\right\}.$$
Since $t(n)=(1-2^{-n})\zeta(n)$, we find that (see also \cite[Theorem 3.3]{Hoffman})
\begin{align}
\Gamma\left(\frac{1-z}{2}\right)=\sqrt{\pi}2^z\exp\left\{\frac{\gamma z}{2}+\sum\limits_{n=2}^\infty \frac{t(n)}{n}z^n\right\}.
\label{Eq:Gamma-2}
\end{align}
Now it is easy to finish the proof.\qed

Similarly, we have a formula for the generating function of the sums of multiple $t$-star values of maximal height.

\begin{cor}
For formal variables $u,w$, we have
\begin{align}
1+\sum\limits_{k\geq 2n,n\geq 1}G_0^{\star}(k,n,n)u^{k-2n}w^n=\exp\left\{\sum\limits_{n=2}^\infty\frac{t(n)}{n}((x^\star)^n+(y^{\star})^n-u^n)\right\},
\label{Eq:MaxHeight-Star}
\end{align}
where $x^{\star}$ and $y^{\star}$ are determined by $x^{\star}+y^{\star}=u$ and $x^{\star}y^{\star}=-w$.
\end{cor}

\proof Setting $v=0$ in Theorem \ref{Thm:G-0-kns-Star}, we get $\alpha^{\star}+\beta^{\star}=3-\frac{u}{2}$ and $\alpha^{\star}\beta^{\star}=\frac{1}{4}(9-3u-w)$. Let $x^{\star}=3-2\alpha^{\star}$ and $y^{\star}=3-2\beta^{\star}$, then $x^{\star}+y^{\star}=u$ and $x^{\star}y^{\star}=-w$. Using the summation formula \eqref{Eq:3F2-1} with $a=\frac{1-u}{2}$, $b=\frac{1}{2}$ and $c=\alpha^{\star}$, we find that
\begin{align*}
&\frac{1}{1-u-w}\,_3F_2\left({\frac{1-u}{2},\frac{1}{2},1\atop \alpha^{\star},\beta^{\star}};1\right)\\
=&\frac{1}{1-u-w}\frac{-(\beta^{\star}-1)}{\left(\alpha^{\star}-\frac{3}{2}\right)\left(\beta^{\star}-\frac{3}{2}\right)}\left(1-\alpha^{\star}+\frac{\Gamma\left(\alpha^{\star}\right)\Gamma\left(\beta^{\star}-1\right)}{\Gamma\left(\frac{1-u}{2}\right)\Gamma\left(\frac{1}{2}\right)}\right)\\
=&-\frac{1}{w}+\frac{1}{w\sqrt{\pi}}\frac{\Gamma\left(\frac{1-x^{\star}}{2}\right)\Gamma\left(\frac{1-y^{\star}}{2}\right)}{\Gamma\left(\frac{1-u}{2}\right)}.
\end{align*}
Now the result follows from \eqref{Eq:Gamma-2}.\qed

From \eqref{Eq:MaxHeight} and \eqref{Eq:MaxHeight-Star}, we see that
$$\left(1+\sum\limits_{k\geq 2n,n\geq 1}G_0(k,n,n)u^{k-2n}w^n\right)\left(1+\sum\limits_{k\geq 2n,n\geq 1}(-1)^nG_0^{\star}(k,n,n)u^{k-2n}w^n\right)=1.$$
Also, setting $u=0$ in \eqref{Eq:MaxHeight} and \eqref{Eq:MaxHeight-Star}, we get
\begin{align*}
&1+\sum\limits_{n=1}^{\infty}t(\{2\}^n)w^n=\exp\left\{\sum\limits_{n=1}^\infty\frac{(-1)^{n-1}}{n}t(2n)w^n\right\},\\
&1+\sum\limits_{n=1}^{\infty}t^{\star}(\{2\}^n)w^n=\exp\left\{\sum\limits_{n=1}^\infty\frac{1}{n}t(2n)w^n\right\}.
\end{align*}
The above formulas  can be deduced from the identities \cite{IKOO,IKZ}
\begin{align*}
&\frac{1}{1-z_kw}=\exp_{\ast}\left(\sum\limits_{n=1}^\infty \frac{(-1)^{n-1}}{n}z_{nk}w^n\right),\\
&\left(\frac{1}{1+z_kt}\right)\ast S\left(\frac{1}{1-z_kt}\right)=1
\end{align*}
in the harmonic shuffle algebra.

\subsection{A weighted sum formula}

Let $I_0(k,n)$ be the set of admissible indices of weight $k$ and depth $n$, and define
$$G_0(k,n)=\sum\limits_{\mathbf{k}\in I_0(k,n)}t(\mathbf{k}),\quad G_0^{\star}(k,n)=\sum\limits_{\mathbf{k}\in I_0(k,n)}t^{\star}(\mathbf{k}),$$
which are the sums of multiple $t$-values and multiple $t$-star values with fixed weight $k$ and depth $n$, respectively. By setting $w=uv$ and then $v=2u$ or $v=-2u$ in Theorem \ref{Thm:G-0-kns} and Theorem \ref{Thm:G-0-kns-Star}, we obtain the following results.

\begin{prop}
For a formal variable $u$, we have
\begin{align}
&\sum\limits_{k=2}^\infty\left(\sum\limits_{n=1}^{k-1}2^{n-1}G_0(k,n)\right)u^{k-2}=\sum\limits_{k=2}^\infty\left(\sum\limits_{n=1}^{k-1}(-1)^{k-n-1}2^{n-1}G_0^{\star}(k,n)\right)u^{k-2}\nonumber\\
=&2^{2u}t(2)\exp\left(\sum\limits_{n=2}^\infty\frac{4(2^{n-1}-1)}{n(2^n-1)}t(n)u^n\right).
\label{Eq:Sum-WeightDepth}
\end{align}
\end{prop}

\proof If $w=uv$ in Theorem \ref{Thm:G-0-kns}, we have
$$\alpha,\beta=\frac{1-u+v}{2},\frac{1}{2}.$$
Hence we get the expression for the generating function of sums of multiple $t$-values with fixed weight and depth:
$$\sum\limits_{k>n\geq 1}G_0(k,n)u^{k-n-1}v^{n-1}=\frac{1}{1-u}\,_3F_2\left({\frac{1-u+v}{2},\frac{1}{2},1\atop \frac{3-u}{2},\frac{3}{2}};1\right).$$
Setting $v=2u$, then we have
$$\sum\limits_{k>n\geq 1}2^{n-1}G_0(k,n)u^{k-2}=\frac{1}{1-u}\,_3F_2\left({\frac{1+u}{2},\frac{1}{2},1\atop \frac{3-u}{2},\frac{3}{2}};1\right).$$
Using Dixon's summation formula \cite[7.4.4.21]{PBA}
\begin{align*}
&\,_3F_2\left({a,b,c\atop 1+a-b,1+a-c};1\right)\\
=&\frac{\sqrt{\pi}}{2^a}\frac{\Gamma(1+a-b)\Gamma(1+a-c)\Gamma\left(1+\frac{a}{2}-b-c\right)}{\Gamma\left(\frac{1+a}{2}\right)\Gamma\left(1+\frac{a}{2}-b\right)\Gamma\left(1+\frac{a}{2}-c\right)\Gamma(1+a-b-c)},\quad \Re(a-2b-2c)>-2,
\end{align*}
with $a=1$, $b=\frac{1}{2}$ and $c=\frac{1+u}{2}$, we get
\begin{align*}
\frac{1}{1-u}\,_3F_2\left({\frac{1+u}{2},\frac{1}{2},1\atop \frac{3-u}{2},\frac{3}{2}};1\right)
=&\frac{1}{1-u}\frac{\sqrt{\pi}}{2}\frac{\Gamma\left(\frac{3}{2}\right)\Gamma\left(\frac{3-u}{2}\right)\Gamma\left(\frac{1-u}{2}\right)}{\Gamma\left(1-\frac{u}{2}\right)\Gamma\left(1-\frac{u}{2}\right)}\\
=&\frac{\pi}{8}\frac{\Gamma\left(\frac{1-u}{2}\right)^2}{\Gamma\left(1-\frac{u}{2}\right)^2}.
\end{align*}
Then using \eqref{Eq:Gamma}, \eqref{Eq:Gamma-2}, the relation $\zeta(n)=(1-2^{-n})^{-1}t(n)$ and the fact $t(2)=\frac{\pi^2}{8}$, we get
$$\sum\limits_{k=2}^\infty\left(\sum\limits_{n=1}^{k-1}2^{n-1}G_0(k,n)\right)u^{k-2}=2^{2u}t(2)\exp\left(\sum\limits_{n=2}^\infty\frac{4(2^{n-1}-1)}{n(2^n-1)}t(n)u^n\right).$$

Similarly, setting $w=uv$ in Theorem \ref{Thm:G-0-kns-Star}, we have
$$\alpha^{\star},\beta^{\star}=\frac{3-u-v}{2},\frac{3}{2}.$$
Hence we find
$$\sum\limits_{k>n\geq 1}G_0^{\star}(k,n)u^{k-n-1}v^{n-1}=\frac{1}{1-u-v}\,_3F_2\left({\frac{1-u}{2},\frac{1}{2},1\atop \frac{3-u-v}{2},\frac{3}{2}};1\right).$$
Let $v=-2u$, then we get
$$\sum\limits_{k>n\geq 1}(-2)^{n-1}G_0^{\star}(k,n)u^{k-2}=\frac{1}{1+u}\,_3F_2\left({\frac{1-u}{2},\frac{1}{2},1\atop \frac{3+u}{2},\frac{3}{2}};1\right).$$
Now it is easy to finish the proof.\qed

Expanding the right-hand side of \eqref{Eq:Sum-WeightDepth}, we get the following weighted sum formula.

\begin{cor}
For any integer $k\geq 2$, we have
\begin{align*}
&\sum\limits_{n=1}^{k-1}2^{n-1}G_0(k,n)=\sum\limits_{n=1}^{k-1}(-1)^{k-n-1}2^{n-1}G_0^{\star}(k,n)\\
=&\sum\limits_{n+n_1+\cdots+n_m=k-2\atop n,m\geq 0,n_1,\ldots,n_m\geq 2}\frac{2^{n+2m}(2^{n_1-1}-1)\cdots(2^{n_m-1}-1)}{n!m!n_1\cdots n_m(2^{n_1}-1)\cdots(2^{n_m}-1)}t(2)t(n_1)\cdots t(n_m)\log^n 2.
\end{align*}
\end{cor}


\section{Proofs of Theorem \ref{Thm:G-0-kns} and Theorem \ref{Thm:G-0-kns-Star}}\label{Sec:Proof}

The proofs of Theorem \ref{Thm:G-0-kns} and Theorem \ref{Thm:G-0-kns-Star} are similar as that of Ohno-Zagier relation for multiple zeta values in \cite{Ohno-Zagier}.

As in \cite{Hoffman}, for an index $\mathbf{k}=(k_1,\ldots,k_n)$, define
\begin{align*}
&\mathcal{L}_{\mathbf{k}}(z)=\mathcal{L}_{k_1,\ldots,k_n}(z)=\sum\limits_{m_1>\cdots>m_n>0\atop m_i:\text{odd}}\frac{z^{m_1}}{m_1^{k_1}\cdots m_n^{k_n}},\\
&\mathcal{L}_{\mathbf{k}}^{\star}(z)=\mathcal{L}_{k_1,\ldots,k_n}^{\star}(z)=\sum\limits_{m_1\geq\cdots\geq m_n>0\atop m_i:\text{odd}}\frac{z^{m_1}}{m_1^{k_1}\cdots m_n^{k_n}}.
\end{align*}
Then $\mathcal{L}_{\mathbf{k}}(z)$ and $\mathcal{L}_{\mathbf{k}}^{\star}(z)$ converge absolutely for $|z|<1$. If $k_1>1$, we have $\mathcal{L}_{\mathbf{k}}(1)=t(\mathbf{k})$ and $\mathcal{L}_{\mathbf{k}}^{\star}(1)=t^{\star}(\mathbf{k})$. From \cite[Lemma 5.1]{Hoffman}, we know that
\begin{align}
\frac{d}{dz}\mathcal{L}_{k_1,\ldots,k_n}(z)=\begin{cases}
\frac{1}{z}\mathcal{L}_{k_1-1,k_2,\ldots,k_n}(z), & \text{if\;} k_1>1,\\
\frac{z}{1-z^2}\mathcal{L}_{k_2,\ldots,k_n}(z), & \text{if\;} n\geq 2, k_1=1,\\
\frac{1}{1-z^2}, & \text{if\;} n=k_1=1.
\end{cases}
\label{Eq:d-MPL}
\end{align}
Similarly, we have
\begin{align}
\frac{d}{dz}\mathcal{L}_{k_1,\ldots,k_n}^{\star}(z)=\begin{cases}
\frac{1}{z}\mathcal{L}_{k_1-1,k_2,\ldots,k_n}^{\star}(z), & \text{if\;} k_1>1,\\
\frac{1}{z(1-z^2)}\mathcal{L}_{k_2,\ldots,k_n}^{\star}(z), & \text{if\;} n\geq 2, k_1=1,\\
\frac{1}{1-z^2}, & \text{if\;} n=k_1=1.
\end{cases}
\label{Eq:d-MPL-Star}
\end{align}
One can also obtain \eqref{Eq:d-MPL} and \eqref{Eq:d-MPL-Star} from the following iterated integral representations:
\begin{align*}
&\mathcal{L}_{k_1,\ldots,k_n}(z)=\int_0^z\left(\frac{dt}{t}\right)^{k_1-1}\frac{tdt}{1-t^2}\cdots\left(\frac{dt}{t}\right)^{k_{n-1}-1}\frac{tdt}{1-t^2}\left(\frac{dt}{t}\right)^{k_n-1}\frac{dt}{1-t^2},\\
&\mathcal{L}_{k_1,\ldots,k_n}^{\star}(z)=\int_0^z\left(\frac{dt}{t}\right)^{k_1-1}\frac{dt}{t(1-t^2)}\cdots\left(\frac{dt}{t}\right)^{k_{n-1}-1}\frac{dt}{t(1-t^2)}\left(\frac{dt}{t}\right)^{k_n-1}\frac{dt}{1-t^2}.
\end{align*}

\subsection{A proof of Theorem \ref{Thm:G-0-kns}}

For nonnegative integers $k,n,s$, denote by $I(k,n,s)$ the set of indices of weight $k$, depth $n$ and height $s$, and define the sums
$$G(k,n,s;z)=\sum\limits_{\mathbf{k}\in I(k,n,s)}\mathcal{L}_{\mathbf{k}}(z),\quad G_0(k,n,s;z)=\sum\limits_{\mathbf{k}\in I_0(k,n,s)}\mathcal{L}_{\mathbf{k}}(z).$$
If the indices set is empty, the sum is treated as zero. And we also set $G(0,0,0;z)=1$.  Note that if $k\geq n+s$ and $n\geq s\geq 1$, we have
$$G_0(k,n,s;1)=G_0(k,n,s).$$
For integers $k,n,s$, using \eqref{Eq:d-MPL}, we have
\begin{description}
  \item[(1)] if $k\geq n+s$ and $n\geq s\geq 1$,
  \begin{align}
  &\frac{d}{dz}G_0(k,n,s;z)\nonumber\\
  =&\frac{1}{z}\left[G(k-1,n,s-1;z)+G_0(k-1,n,s;z)-G_0(k-1,n,s-1;z)\right],
  \label{Eq:d-G-0}
  \end{align}
  \item[(2)] if $k\geq n+s$, $n\geq s\geq 0$ and $n\geq 2$,
  \begin{align}
  \frac{d}{dz}\left[G(k,n,s;z)-G_0(k,n,s;z)\right]=\frac{z}{1-z^2}G(k-1,n-1,s;z).
  \label{Eq:d-G-G0}
  \end{align}
\end{description}

Then we define the generating functions
\begin{align*}
\Phi(z)=&\Phi(u,v,w;z)=\sum\limits_{k,n,s\geq 0}G(k,n,s;z)u^{k-n-s}v^{n-s}w^{s},\\
=&1+\mathcal{L}_1(z)v+\sum\limits_{k\geq 2}\mathcal{L}_k(z)u^{k-2}w+\sum\limits_{k\geq n+s\atop n\geq s\geq 0,n\geq 2}G(k,n,s;z)u^{k-n-s}v^{n-s}w^{s},\\
\Phi_0(z)=&\Phi_0(u,v,w;z)=\sum\limits_{k,n,s\geq 0}G_0(k,n,s;z)u^{k-n-s}v^{n-s}w^{s-1}\\
=&\sum\limits_{k\geq n+s\atop n\geq s\geq 1}G_0(k,n,s;z)u^{k-n-s}v^{n-s}w^{s-1}.
\end{align*}
Using \eqref{Eq:d-MPL}, \eqref{Eq:d-G-0} and \eqref{Eq:d-G-G0}, we see that
\begin{align*}
&\frac{d}{dz}\Phi_0(z)=\frac{1}{vz}(\Phi(z)-1-w\Phi_0(z))+\frac{u}{z}\Phi_0(z),\\
&\frac{d}{dz}(\Phi(z)-w\Phi_0(z))=\frac{vz}{1-z^2}\Phi(z)+\frac{v}{1+z}.
\end{align*}
Eliminating $\Phi(z)$, we obtain the differential equation satisfied by $\Phi_0(z)$.

\begin{prop}
$\Phi_0=\Phi_0(z)$ satisfies the following differential equation
\begin{align}
z(1-z^2)\Phi_0''+[(1-u)(1-z^2)-vz^2]\Phi_0'+(uv-w)z\Phi_0=1.
\label{Eq:G0-DE}
\end{align}
\end{prop}

We want to find the unique power series solution $\Phi_0(z)=\sum\limits_{n=1}^\infty a_nz^n$. From \eqref{Eq:G0-DE}, we see that $a_1=\frac{1}{1-u}$, $a_2=0$ and
$$a_{n+1}=\frac{(n-1)(n-2)+(1-u+v)(n-1)-(uv-w)}{(n+1-u)(n+1)}a_{n-1},\quad n\geq 2.$$
Hence for any $n\geq 1$, we have $a_{2n}=0$ and
$$a_{2n+1}=\frac{\left(n-\frac{1}{2}\right)(n-1)+\left(\frac{1}{2}-\frac{u}{2}+\frac{v}{2}\right)\left(n-\frac{1}{2}\right)-\frac{1}{4}(uv-w)}{\left(n+\frac{1}{2}-\frac{u}{2}\right)\left(n+\frac{1}{2}\right)}a_{2n-1}.$$
Since $(\alpha-1)+(\beta-1)=-1-\frac{u}{2}+\frac{v}{2}$ and $(\alpha-1)(\beta-1)=\frac{1}{4}(1+u-v-uv+w)$, we have
$$a_{2n+1}=\frac{(n+\alpha-1)(n+\beta-1)}{\left(n+\frac{1}{2}-\frac{u}{2}\right)\left(n+\frac{1}{2}\right)}a_{2n-1}=\frac{(\alpha)_n(\beta)_n}{\left(\frac{3-u}{2}\right)_n\left(\frac{3}{2}\right)_n}\frac{1}{1-u}.$$
Therefore, we can represent $\Phi_0(z)$ by the generalized hypergeometric function $\,_3F_2$ as displaying in the following theorem.

\begin{thm}\label{Thm:G-0-z-kns}
We have
$$\Phi_0(k,n,s;z)=\frac{z}{1-u}\,_3F_2\left({\alpha,\beta,1\atop \frac{3-u}{2},\frac{3}{2}};z^2\right).$$
\end{thm}

Finally, setting $z=1$ in Theorem \ref{Thm:G-0-z-kns}, we get Theorem \ref{Thm:G-0-kns}. \qed

\subsection{A proof of Theorem \ref{Thm:G-0-kns-Star}}

Similarly, for nonnegative integers $k,n,s$, we define the sums
$$G^{\star}(k,n,s;z)=\sum\limits_{\mathbf{k}\in I(k,n,s)}\mathcal{L}_{\mathbf{k}}^{\star}(z),\quad G_0^{\star}(k,n,s;z)=\sum\limits_{\mathbf{k}\in I_0(k,n,s)}\mathcal{L}_{\mathbf{k}}^{\star}(z)$$
with $G^{\star}(0,0,0;z)=1$. Using \eqref{Eq:d-MPL-Star}, we have
\begin{description}
  \item[(1)] if $k\geq n+s$ and $n\geq s\geq 1$,
  \begin{align}
  &\frac{d}{dz}G_0^{\star}(k,n,s;z)\nonumber\\
  =&\frac{1}{z}\left[G^{\star}(k-1,n,s-1;z)+G_0^{\star}(k-1,n,s;z)-G_0^{\star}(k-1,n,s-1;z)\right],
  \label{Eq:d-G-0-Star}
  \end{align}
  \item[(2)] if $k\geq n+s$, $n\geq s\geq 0$ and $n\geq 2$,
  \begin{align}
  \frac{d}{dz}\left[G^{\star}(k,n,s;z)-G_0^{\star}(k,n,s;z)\right]=\frac{1}{z(1-z^2)}G^{\star}(k-1,n-1,s;z).
  \label{Eq:d-G-G0-Star}
  \end{align}
\end{description}
We define the generating functions
\begin{align*}
\Phi^{\star}(z)=&\Phi^{\star}(u,v,w;z)=\sum\limits_{k,n,s\geq 0}G^{\star}(k,n,s;z)u^{k-n-s}v^{n-s}w^{s},\\
\Phi_0^{\star}(z)=&\Phi_0^{\star}(u,v,w;z)=\sum\limits_{k,n,s\geq 0}G_0^{\star}(k,n,s;z)u^{k-n-s}v^{n-s}w^{s-1}.
\end{align*}
Then using \eqref{Eq:d-MPL-Star}, \eqref{Eq:d-G-0-Star} and \eqref{Eq:d-G-G0-Star}, we see that
\begin{align*}
&\frac{d}{dz}\Phi_0^{\star}(z)=\frac{1}{vz}(\Phi^{\star}(z)-1-w\Phi_0^{\star}(z))+\frac{u}{z}\Phi_0^{\star}(z),\\
&\frac{d}{dz}(\Phi^{\star}(z)-w\Phi_0^{\star}(z))=\frac{v}{z(1-z^2)}\Phi^{\star}(z)-\frac{v}{z(1+z)}.
\end{align*}
Eliminating $\Phi^{\star}(z)$, we get the differential equation satisfied by $\Phi_0^{\star}(z)$.

\begin{prop}
$\Phi_0^{\star}=\Phi_0^{\star}(z)$ satisfies the following differential equation
\begin{align}
z^2(1-z^2)(\Phi_0^{\star})''+[(1-u)z(1-z^2)-vz](\Phi_0^{\star})'+(uv-w)\Phi_0^{\star}=z.
\label{Eq:G0-DE-Star}
\end{align}
\end{prop}

Assume that $\Phi_0^{\star}(z)=\sum\limits_{n=1}^\infty a_n^{\star}z^n$. Using \eqref{Eq:G0-DE-Star}, we find that $a_1^{\star}=\frac{1}{1-u-v+uv-w}$, $a_2^{\star}=0$ and
$$a_{n}^{\star}=\frac{(n-2)(n-2-u)}{n(n-1)+n(1-u-v)+uv-w}a_{n-2}^{\star},\quad n\geq 3.$$
Hence for any $n\geq 1$, we have $a_{2n}^{\star}=0$ and
$$a_{2n+1}^{\star}=\frac{\left(n-\frac{1}{2}-\frac{u}{2}\right)\left(n-\frac{1}{2}\right)}{(n+\alpha^{\star}-1)(n+\beta^{\star}-1)}a_{2n-1}^{\star}=\frac{\left(\frac{1-u}{2}\right)_n\left(\frac{1}{2}\right)_n}{(\alpha^{\star})_n(\beta^{\star})_n}
\frac{1}{1-u-v+uv-w}.$$
Therefore, we have the following theorem.

\begin{thm}\label{Thm:G-0-z-kns-Star}
We have
$$\Phi_0^{\star}(k,n,s;z)=\frac{z}{1-u-v+uv-w}\,_3F_2\left({\frac{1-u}{2},\frac{1}{2},1\atop \alpha^{\star},\beta^{\star}};z^2\right).$$
\end{thm}

Finally, setting $z=1$ in Theorem \ref{Thm:G-0-z-kns-Star}, we get Theorem \ref{Thm:G-0-kns-Star}.\qed

\end{document}